\documentclass[titlepage,11pt]{article}
\usepackage{latexsym}
\usepackage{amsfonts}
\usepackage{amsmath}
\usepackage{mathtools}
\usepackage{tikz}
\usetikzlibrary{shapes.geometric}
\newcommand\blackslug{\hbox{\hskip 1pt \vrule width 4pt height 8pt depth 1.5pt
		\hskip 1pt}}
\newcommand\bbox{\hfill \quad \blackslug \bigbreak}
\newcommand{\erh}{Erd\H{o}s-Hajnal}
\def\DD{\hbox{-}}

\def\LL{,\ldots,}
\newcommand{\vare}{\varepsilon}

\newcommand{\cupcup}{\cup \cdots\cup}

\def\poly{\operatorname{poly}}
\newcommand{\mac}{\mathcal}
\def\ind{\operatorname{ind}}

\DeclarePairedDelimiter\ceil{\lceil}{\rceil}%

\DeclarePairedDelimiter\abs{\lvert}{\rvert}%

\title{Induced subgraph density. III. Cycles and subdivisions}
\author{
	Tung Nguyen\thanks{Supported by AFOSR grants
		A9550-19-1-0187 and FA9550-22-1-0234, and by NSF grants  DMS-1800053 and DMS-2154169.}\\
	Princeton University,\\ Princeton, NJ 08544, USA
	\and
	Alex Scott\thanks{Supported by EPSRC grant EP/X013642/1}\\
	University of Oxford, \\
	Oxford, UK
	\and
	Paul Seymour\thanks{Supported by AFOSR grants
		A9550-19-1-0187 and FA9550-22-1-0234, and by NSF grants  DMS-1800053 and DMS-2154169.}\\
	Princeton University,\\ Princeton, NJ 08544, USA}

\date{July 4, 2022; revised \today}

\newtheorem{thm}{}[section]

\newcommand{\Proof}{\noindent{\bf Proof.}\ \ }

\begin{document}
	\maketitle
	\begin{abstract}
		We show that for every two cycles $C,D$, there exists $c>0$ such that if $G$ is both $C$-free and $\overline{D}$-free then
		$G$ has a clique or stable set of size at least $|G|^c$. (``$H$-free''            
		means 
		with no induced subgraph isomorphic to $H$, and $\overline{D}$ denotes the complement graph of $D$.)
		Since the five-vertex cycle $C_5$ is isomorphic to its complement, this extends the earlier result that $C_5$ satisfies the Erd\H{o}s-Hajnal 
		conjecture.
		It also unifies and strengthens several other results.
		
		The results for cycles are special cases of results for subdivisions, as follows. 
		Let $H,J$ be obtained from smaller graphs by subdividing every edge exactly twice.
		We will prove that there exists $c>0$ such that if $G$ is both $H$-free and $\overline{J}$-free then
		$G$ has a clique or stable set of size at least $|G|^c$. And the same holds if $H$ and/or $J$ is obtained from a graph by 
		choosing a forest $F$ and subdividing every edge not in $F$ at least five times.  Our proof uses the framework of iterative sparsification developed in other papers of this series.
		
		Along the way, we will also give a short and simple proof strengthening a celebrated result of Fox and Sudakov, that says that for all $H$,
		every $H$-free graph contains either a large stable set or a large complete bipartite subgraph.

	\end{abstract}
	
	\section{Introduction}
	
	A graph is {\em $H$-free} if it has no induced subgraph isomorphic to $H$; and if
	$\mathcal{H}$  is a set of graphs, then a graph $G$ is {\em $\mathcal{H}$-free} if it is $H$-free for all $H\in \mathcal{H}$.
	We say a set $\mathcal{H}$ has the {\em Erd\H{o}s-Hajnal property} if there exists $c>0$
	such that for every $\mathcal{H}$-free graph $G$, there is a clique or stable set of $G$ with 
	cardinality at least $|G|^c$; and a graph $H$ has the property if $\{H\}$ does.
	
	A famous old conjecture of Erd\H{o}s and Hajnal~\cite{EH0, EH} says that every graph has the 
	Erd\H{o}s-Hajnal property; or equivalently, every non-null set of graphs has the property.
	Despite extensive efforts, the conjecture has only been proved for a small family of graphs $H$. By a theorem of Alon, Pach and 
	Solymosi~\cite{APS}, graphs that are made by vertex-substitution from other graphs with the property also have the property,
	so we could confine our attention to {\em prime} graphs, those that cannot be built from smaller graphs by vertex-substitution. 
	But until recently, there were only three  prime graphs with more than two vertices that were known to have the property: the 
	four-vertex path, the bull~\cite{safra} and the five-vertex cycle~\cite{C5}. (Very recently, more graphs have been added to this list, including
	the five-vertex path~\cite{density7} and infinitely many other prime graphs~\cite{density4}; see also \cite{density5, density6,equiv} for other recent progress.)
	
	There has been significant progress when more than one graph is excluded
	(we review this in detail in the next section).  For example, it was shown in~\cite{C5} that $\{H,\overline H\}$ has the Erd\H{o}s-Hajnal property if $H$ is a cycle of length at most seven~\cite{C5}, and extending this to 
	longer cycles was mentioned as a nice open question.  Here we will prove a substantially stronger result:
	
	\begin{thm}\label{holeandantihole}
		Let $H$ and $J$ be cycles.  Then 
		$\{H,\overline{J}\}$ has the Erd\H{o}s-Hajnal property.
	\end{thm}
	
	We will also show that:
	
	\begin{thm}\label{2subdivisionsimple}
		Let $H$ be obtained from the complete graph $K_n$ by subdividing every edge twice. Then $\{H,\overline{H}\}$ has the
		Erd\H{o}s-Hajnal property.
	\end{thm}
	
	Both these results are consequence of more general theorems, which we give below.  We discuss related work in section \ref{strongehp} and then give our results in section \ref{ourresults}.

	\section{Related work}\label{strongehp}
	
	There has been a substantial body of work proving the Erd\H{o}s-Hajnal property when two or more graphs are excluded.  Some of this has proceeded by proving a stronger property: 
	we say that a set of graphs
	$\mathcal{H}$ has the {\em strong Erd\H{o}s-Hajnal property} if there exists $c>0$
	such that for every $\mathcal{H}$-free graph $G$ with at least two vertices, there are disjoint sets $A,B\subseteq V(G)$ with size at least $c|G|$ such that the pair
	$(A,B)$ is either complete or anticomplete (i.e.~either all possible edges between $A$ and $B$ are present in $G$, or there are no edges between $A$ and $B$).
	It is not hard to show that the strong Erd\H{o}s-Hajnal property implies the Erd\H{o}s-Hajnal property.
	
	Which families $\mathcal H$ satisfy the strong Erd\H{o}s-Hajnal property?  By considering sparse random graphs, it follows that every finite set $\mathcal H$ with the strong Erd\H{o}s-Hajnal property must contain both a forest and the complement of a forest.  It was
	shown in~\cite{pure1} that every such $\mathcal H$ has the strong Erd\H os-Hajnal property:
	
	\begin{thm}\label{ffbar}
		Let $H$ and $J$ be forests.  Then
		$\{H,\overline{J}\}$ has the strong
		Erd\H{o}s-Hajnal property.
	\end{thm}
	
	This characterizes finite sets $\mathcal H$ with the strong Erd\H{o}s-Hajnal property, and shows that excluding both a forest and the complement of a forest implies the Erd\H{o}s-Hajnal property.
	(By contrast, much less is known if we forbid just a forest $H$: for example, the Erd\H{o}s-Hajnal property was only recently proved 
	in the case $H=P_5$~\cite{density7}, and only the ``near Erd\H{o}s-Hajnal'' property is known if $H$ is a longer path~\cite{density5}.)
	
	If $\mathcal H$ has the strong Erd\H{o}s-Hajnal property but does not include both a forest and the complement of a forest, then $\mathcal H$ must be infinite.  There has been progress here as well.
	Bonamy, Bousquet and Thomass\'e~\cite{bonamy} showed that excluding all long holes and antiholes gives the strong Erd\H{o}s-Hajnal property:
	\begin{thm}\label{bonamy}
		Let $k\ge 3$, and let $\mac H$ contain all cycles of length at least $k$ and their complements. Then $\mac H$ has the strong 
		Erd\H{o}s-Hajnal property.
	\end{thm}
	It is interesting to compare this with \ref{holeandantihole}.
	One can show it is necessary to exclude infinitely many cycles and infinitely many complements of cycles to obtain 
	the strong Erd\H{o}s-Hajnal property.  By contrast,
	\ref{holeandantihole} shows that excluding a {\em single} cycle and a {\em single} cycle complement is enough to give the Erd\H{o}s-Hajnal property.
	
	A strengthening of \ref{bonamy} was also known, but first we need some definitions.
	{\em Subdividing} an edge $uv$ means deleting the edge, and adding a path between $u,v$ whose internal vertices are new. If the path
	has length $k+1$ this is called {\em $k$-subdividing} the edge; and if the path has length at least $k+1$ it is called
	{\em $(\ge k)$-subdividing}. A graph obtained from another graph $H$ by subdividing some of the edges of $H$ is called a {\em subdivision}
	of $H$, and it is a {\em $k$-subdivision} or {\em $(\ge k)$-subdivision} if every edge is $k$-subdivided or every edge is $(\ge k)$-subdivided, 
	respectively. 
	
	The following substantial strengthening of \ref{bonamy}  was proved in~\cite{pure2}:
	\begin{thm}\label{pure2}
		Let $H$ be a graph, and let $\mathcal{H}$ contain all subdivisions of $H$ and their complements. Then $\mac H$ has the strong
		Erd\H{o}s-Hajnal property.
	\end{thm}
	Once again, we need to exclude infinitely many induced subdivisions of $H$, and the complements of infinitely many induced subdivisions.  By contrast,~\ref{2subdivisionsimple} shows that the Erd\H{o}s-Hajnal property holds if we exclude a {\em single} subdivision and its complement.  In fact, we prove a more general result (\ref{swiss}), which we discuss in the next section.
	
	\section{Results}\label{ourresults}
	
	\ref{holeandantihole} and~\ref{2subdivisionsimple} are both consequences of the following more general results.
	\begin{thm}\label{2subdivision}
		Let $H, J$ be 2-subdivisions of multigraphs $H_0,J_0$ respectively. Then $\{H,\overline{J}\}$ has the
		Erd\H{o}s-Hajnal property.
	\end{thm}
	\begin{thm}\label{5subdivision}
		Let $H$ be obtained from a multigraph $H_0$ by choosing a forest $F$ of $H_0$, and $(\ge 5)$-subdividing every edge 
		of $H_0$ not in $E(F)$. Let $J$ be constructed similarly. Then $\{H,\overline{J}\}$ has the
		Erd\H{o}s-Hajnal property.
	\end{thm}
	Let us mention that we can also prove a variant of \ref{5subdivision} (this will appear in Tung Nguyen's thesis~\cite{tungthesis}):
	\begin{thm}\label{4subdivision}
		Let $H, J$ be $(\ge 4)$-subdivisions of multigraphs $H_0, J_0$ respectively.
		Then $\{H,\overline{J}\}$ has the
		Erd\H{o}s-Hajnal property.
	\end{thm}
	We do not give its proof here, and the proof method is different from those used in this paper.

	There is a common generalization of \ref{2subdivision} and \ref{5subdivision}, which is our objective. 
	Let $F$ be a forest, and let $s,t\ge 1$ be integers. The graph $F^s_t$ is obtained as follows. Let us:
	\begin{itemize}
		\item add a set $X$
		of $t$ new vertices to $F$;
		\item for all distinct $u,v\in X$, add $s$ parallel edges with ends $u,v$,
		and for each $u\in V(F)$ and $v\in X$, add $s$ parallel edges with ends $u,v$,
		making a multigraph;
		\item subdivide exactly twice every edge of this multigraph with an end in $X$ (that is, all its edges except those of $F$).
	\end{itemize}
	If $H$ is a graph such that $H=F^s_t$ for some such $F,s,t$, we call $H$ a {\em Swiss Army graph}.
	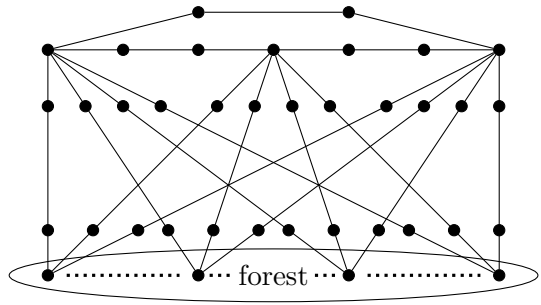
\begin{figure}[h!]
		\centering
		
		\begin{tikzpicture}[scale=1/2,auto=left]
			
			\draw (0,0) ellipse (200 pt and 20 pt);
			\node at (0,0) {forest};
			\tikzstyle{every node}=[inner sep=1.5pt, fill=black,circle,draw]
			\node (f1) at (-6,0) {};
			\node (f2) at (-2,0) {};
			\node (f3) at (2,0) {};
			\node (f4) at (6,0) {};
			\node (a1) at (-6,6) {};
			\node (a2) at (0,6) {};
			\node (a3) at (6,6) {};
			
			\draw (f1)--(a1);
			\draw (f1)--(a2);
			\draw (f1)--(a3);
			\draw (f2)--(a1);
			\draw (f2)--(a2);
			\draw (f2)--(a3);
			\draw (f3)--(a1);
			\draw (f3)--(a2);
			\draw (f3)--(a3);
			\draw (f4)--(a1);
			\draw (f4)--(a2);
			\draw (f4)--(a3);
			
			\node (b11) at (-6,4.5) {};
			\node (b12) at (-1.5,4.5) {};
			\node (b13) at (3,4.5) {};
			\node (b21) at (-5,4.5) {};
			\node (b22) at (-.5,4.5) {};
			\node (b23) at (4,4.5) {};
			\node (b31) at (-4,4.5) {};
			\node (b32) at (.5,4.5) {};
			\node (b33) at (5,4.5) {};
			\node (b41) at (-3,4.5) {};
			\node (b42) at (1.5,4.5) {};
			\node (b43) at (6,4.5) {};
			
			\node (c11) at (-6,1.2) {};
			\node (c12) at (-4.8,1.2) {};
			\node (c13) at (-3.6,1.2) {};
			\node (c21) at (-2.8,1.2) {};
			\node (c22) at (-1.6,1.2) {};
			\node (c23) at (-.4,1.2) {};
			\node (c31) at (.4,1.2) {};
			\node (c32) at (1.6,1.2) {};
			\node (c33) at (2.8,1.2) {};
			\node (c41) at (3.6,1.2) {};
			\node (c42) at (4.8,1.2) {};
			\node (c43) at (6,1.2) {};
			
			\node (d12) at (-4,6) {};
			\node (e12) at (-2,6) {};
			\node (d23) at (2,6) {};
			\node (e23) at (4,6) {};
			\node (d13) at (-2,7) {};
			\node (e13) at (2,7) {};
			
			\draw (a1) -- (a3);
			\draw (a1)--(d13)--(e13)--(a3);
			\draw[dotted, very thick] (-5.5,0)--(-2.5,0);
			\draw[dotted, very thick] (2.5,0)--(5.5,0);
			\draw[dotted, very thick] (-1.7,0)--(-1.1,0);
			\draw[dotted, very thick] (1.1,0)--(1.7,0);
			
		\end{tikzpicture}
		
		\caption{A Swiss Army graph with $s=1$ and $t=3$.} \label{fig:swissarmy}
	\end{figure}
	
	We will prove:
	\begin{thm}\label{swiss}
		If $H,J$ are Swiss Army graphs, then $\{H,\overline{J}\}$ has the Erd\H{o}s-Hajnal property.
	\end{thm}
	
	Swiss Army graphs are cumbersome objects, but they contain several useful and interesting features. For instance:
	\begin{itemize}
		\item Every cycle of length at least six is an induced subgraph of a Swiss Army graph. (If the cycle has length at least seven,
		make the forest a path of the right length, while if the cycle has length six take $s=2$.)
		\item Every 2-subdivision of a multigraph
		is an induced subgraph of a Swiss Army graph. (Let the forest be the null graph, and $s$ the
		edge-multiplicity of the multigraph.)
		\item Let $H$ be obtained from
		a multigraph $H_0$ by choosing some forest $F$ of $H_0$, and $(\ge 5)$-subdividing every edge of $H_0$ not in $F$.
		Then $H$ is an induced subgraph of a Swiss Army graph. (By adding leaves to $F$, we may assume every edge of $H_0$
		not in $F$ needs to be subdivided {\em exactly} five times, and then the result is clear: take $t$ to be the number of edges 
		of $H_0$ not in $F$.)
	\end{itemize}
	
	Since we already know the Erd\H{o}s-Hajnal property for $H$-free graphs when $H$ is a cycle of length five or less, the first bullet 
	above shows that  \ref{swiss} implies \ref{holeandantihole}. The second bullet shows that it implies \ref{2subdivision}, and the 
	third that it implies
	\ref{5subdivision}. 
	So now we need to prove \ref{swiss}, and that occupies the remainder of the paper.

	James Davies (private communication) pointed out a nice application of our results, for ``pivot-minors'' (for definitions, see~\cite{davies}). For every graph $J$, there
	is a graph $H$ such that if $G$ does not contain $J$ as a pivot-minor, then it contains neither the 2-subdivision of $H$ nor 
	its complement as an induced subgraph (see Lemmas 2.1 and
	2.3 of~\cite{davies}). Consequently, \ref{2subdivision} implies that for every graph $J$, the class of graphs not containing 
	$J$ as a pivot-minor has the \erh{} property, a result recently proved by Davies in~\cite{davies}. 
	(He actually proved more, the strong \erh{} property.)
	
	We use standard terminology.   Graphs are assumed to be finite, and to have no parallel edges or loops. 
	Occasionally we need to allow parallel edges (but not loops), and in that case we call them {\em multigraphs}; $G[X]$ denotes the
	induced subgraph with vertex set $X$
	of a graph $G$; $|G|$ denotes the number of vertices of $G$; and $\overline{G}$ is the complement graph of $G$.

	%%%%%%%%%%%%%%%%%%%%%%%%%%%%%%%%%%%%%%%%%%%%%%%%%%%%%%%%%%%%%%%%%%%%%%%%%%%%%%%%%%%
	\section{A sketch of the proof}
	
	Let us give some idea of how the proof will work. (In the proof proper, there are numerous constant factors that we will ignore here.)
	We have two Swiss Army graphs $H,J$, and we can assume that $H=J$ without loss of generality. 
	We need to show that 
	every $\{H,\overline{H}\}$-free graph has a stable set or clique of size $|G|^c$, where $c>0$ is some small constant depending on 
	$H$ but not on $G$. It is just as good to show that $\alpha(G)\omega(G)\ge |G|^c$ for all 
	such graphs $G$. Choose $c>0$ very small,
	and suppose it does not work;
	then there is a minimal counterexample $G$, that is, $G$ is ``$c$-critical''. 
	
	We use the strategy of iterative sparsification developed in other papers of this series (see~\cite{density4, density5, density6, density7}).
	Since $G$ is $H$-free, a theorem of R\"odl implies that there is a subset $S\subseteq V(G)$ with size linear in $|G|$, inducing a subgraph
	that is either very sparse or very dense; and by replacing $G$ by its complement we may assume it is very sparse. Consequently, for 
	some 
	$0<y<1$ (at most any positive constant we wish),
	there is a subset $S\subseteq V(G)$ with density at most $y$ and size at least $y^a|G|$ (where $a$ is some large positive constant that we choose for convenience). 
	Say such a value of $y$ is ``good'', and choose a good value of $y$ such that $y^2$ is not good. There must be such a value, because
	there is no good value that is between $|G|^{-c}$ and $|G|^{-2c}$ (because if there were, we could find a large stable set contradicting
	that $G$ is $c$-critical). From now on we just work inside the set $S$.
	
	Now the proof breaks into two parts. We need to find in $G[S]$ an appropriate object, a ``blown-up'' relative of a Swiss Army graph that we
	(temporarily) call a ``template'', that can be described as follows. The graph $H$ equals $F^s_t$ for some $F,s,t$;
	let $F'$ be a forest with no edges and about $y^{-1/16}$ vertices, and consider the graph $(F')^s_t$. 
	Now blow up each vertex of $(F')^s_t$ that belongs to $F'$ into a subset, a ``block'',
	pairwise disjoint, and pairwise sparse  (with sparsity at most $y^{1/6}$ to each other). We want each of the blocks
	to be large (at
	least $\poly(y)|G|$ for some fixed polynomial), but there is no condition on the subgraph induced on each block. 
	That is what we mean by a template.
	
	The first half of the proof is to show that we can find a template in $G[S]$, using the minimality of $y$. 
	Then, once we have such an object, the second half of the proof is to find a copy of $F$ within the union of 
	the blocks, with at most one 
	vertex in each block (that is, a ``rainbow'' copy of $F$). If we can do that, then $F$ can be extended within the template to make a copy of $H$ in $G$, a contradiction
	(which would prove that there is no $\{H,\overline{H}\}$-free $c$-critical graph $G$ after all).

	The first half of the proof, finding the template within $G[S]$, uses a refinement of a result of Fox and Sudakov,
	that we prove in the next 
	section and use in sections \ref{sec:sparseseq} and \ref{sec:blockades}. The theorem says that for any $H$, in any $H$-free graph 
	we can find either a large sparse subset or two large sets of 
	vertices
	that are very dense to each other (and ``large'', ``sparse'' and ``very dense'' can be tuned). We will apply this inside the graph 
	$G[S]$, which is already $y$-sparse, and the first outcome is impossible because of the minimality of $y$. If the second happens, with
	two large sets $A_1,B_1$ say, there are only a few vertices in $S\setminus (A_1\cup B_1)$ that have many neighbours in $A_1$ or in $B_1$ (because $G[S]$ is 
	$y$-sparse), so we can put them aside and repeat. We get a sequence of about $y^{-1/16}$ pairs $A_i,B_i$  of large sets, each pair very dense to each other, and the sets in
	different pairs sparse to each other. Each $B_i$ includes a large stable set (since $G$ is $c$-critical). If we can get a large 
	induced matching within the union of these stable sets, with not many edges joining the same two sets, then we have the template; 
	and otherwise
	there is a large stable set that contradicts the $c$-criticality of $G$.

	The second half uses 
	a result
	of~\cite{pure1}, that says that for every forest $F$, in every sparse $F$-free graph $G$ there are two sets of vertices of linear size
	that are anticomplete (that is, there are no edges between them). We need to replace the $F$-free hypothesis with the weaker hypothesis that there is no rainbow copy of $F$; 
	and to replace the ``two anticomplete sets of linear size'' outcome
	with a ``$\poly(1/y)$ anticomplete sets of $\poly(y)|G|$ size'' outcome. Then we apply this result to the template.
	The ``$\poly(1/y)$ anticomplete sets of $\poly(y)|G|$ size'' outcome would contradict the $c$-criticality of $G$, so there must be 
	a rainbow copy of $F$, which is what we want. This is all done in section \ref{sec:sparsetoanti}.

	%%%%%%%%%%%%%%%%%%%%%%%%%%%%%%%%%%%%%%%%%%%%%%%%%%%%%%%%%%%%%%%%%%%%%%%%%%%%%%%%%%%%
	
	\section{Stable sets and complete bipartite subgraphs}
	
	We need to use a strengthening of a celebrated result of Fox and Sudakov~\cite{foxsudakov}. Their proof was complicated, and used dependent random choice,
	and we begin with giving a simple proof of their theorem. Then we will modify it to prove the strengthening that we need.
	
	The result is 
	an asymmetric weakening of the Erd\H{o}s-Hajnal conjecture in which the
	clique is replaced by a complete bipartite subgraph (not necessarily induced).
	\begin{thm}\label{foxsudakov}{\rm (Fox and Sudakov~\cite{foxsudakov})}
		For every graph $H$, there exists $c>0$ such that
		every $H$-free graph $G$ contains either
		a stable set of size  at least $|G|^c$, or
		a complete bipartite subgraph with both parts of size at least $|G|^c$.
	\end{thm}
	We will give a simple proof of  the following stronger result:
	
	\begin{thm}\label{foxsud1}
		For every graph $H$, and every $\delta>0$, there exists $c>0$ such that the following holds:
		every sufficiently large $H$-free graph $G$ contains either a stable set of size $|G|^c$, or a complete bipartite subgraph with
		parts of cardinality at least $|G|^{1-\delta}$ and $|G|^c$ respectively.
	\end{thm}
	The same method, with a little more effort, also yields a second strengthening of \ref{foxsudakov}, which is what we actually need, 
	but we will prove that separately. If $H,G$ are graphs, a {\em copy} of $H$ in $G$ means an isomorphism from $H$ to an induced subgraph of $G$, and $\ind_H(G)$
	denotes the number of copies of $H$ in $G$.

	The proof of \ref{foxsud1} relies on the following key lemma, which was proved (with different constants) by
	Fox and Sudakov~\cite{foxsudakov}; as they showed, \ref{foxsudakov} follows in a few lines.
	\begin{thm}\label{foxsudmain}
		Let $H$ be a graph, and let $0<\vare<1$. Let $G$ be an $H$-free graph, and let $t>0$ be an integer, with $|G|\ge t$.
		Then either
		\begin{itemize}
			\item there is a stable set of $G$ with size $t$; or
			\item there are disjoint subsets $W_1,W_2$ of $V(G)$, with $|W_1|,|W_2|\ge (2t)^{-|H|^2}\vare^{|H|^2/2}|G|$,
			such that every vertex in $W_1$ is nonadjacent to at most $2\vare|W_2|$ vertices  in $W_2$.
		\end{itemize}
	\end{thm}
	\Proof We may assume that $G$ has no stable set of size $t$, and so it has an edge; and therefore we may assume that $|G|>t^{|H|^2}$, since otherwise the second
	bullet holds (taking $|W_1|=|W_2|=1$).
	So we may assume that $t,|H|\ge 2$.
	Let $k:=|H|$ and $n:=|G|$.
	
	Let $H_0,H_1\LL H_m$ be a sequence of graphs, each with vertex set $V(H)$, where $m=k(k-1)/2$, such that
	for $1\le i\le m$, $H_i$ is obtained
	from $H_{i-1}$ by adding an edge joining two nonadjacent vertices of $H_{i-1}$ (and consequently $H_i$ has $i$ edges),
	and such that one of $H_0\LL H_m$ equals $H$.
	The proof starts with $H_m$ and works downwards.  For each $i\ge1$, we will
	show that if $G$ contains many copies of $H_i$ then either it contains many
	copies of $H_{i-1}$ or we can find a very dense bipartite subgraph.  Since,
	as we shall see, we have many copies of $H_m$ and no copy of $H$, the bipartite outcome must
	occur at some point. 
	We begin with:
	\\
	\\
	(1) {\em  $\ind_{H_m}(G)\ge t^{-k^2}n^k$.}
	\\
	\\
	Every graph with at least $t^k$ vertices has either a stable set of size $t$ or a $k$-clique (that is, a clique of size $k$),
	by one of the standard forms of Ramsey's theorem. By our assumption, $G$ has no stable set of size $t$, and so
	every subset $X\subseteq V(G)$ with $|X|=s$ includes a $k$-clique, where $s=t^k$. Since each $k$-clique is included in only
	$\binom{n-k}{s-k}$ subsets of size $s$, and there are $\binom{n}{s}$ subsets of size $s$ altogether (because
	$n\ge t^k$), it follows that there are at least
	$$\binom{n}{s}/\binom{n-k}{s-k}= \frac{n(n-1)\cdots(n-k+1)}{s(s-1)\cdots(s-k+1)}\ge \left(\frac{n}{s}\right)^k=t^{-k^2}n^k$$
	$k$-cliques in $G$. This proves (1).
	
	\bigskip
	For $0\le i\le m$, let $f(i)=\ind_{H_i}(G)/n^k$. Thus we have seen that
	$f(m)\ge t^{-k^2}$.
	\\
	\\
	(2) {\em For $1\le i\le m$, either $f(i-1)\ge (\vare/4) f(i)$, or there are disjoint subsets $W_1,W_2$ of $V(G)$
		with $|W_1|\cdot|W_2|\ge (f(i)/4)n^2$, such that there are fewer than $\vare|W_1|\cdot |W_2|$ nonedges between $W_1,W_2$.}
	\\
	\\
	We suppose the second outcome is false.
	Let $p,q\in V(H)$ be distinct, such that $H_{i}$ is obtained from $H_{i-1}$ by adding the edge $pq$.
	Let $J=H_{i}\setminus \{p,q\}$, and let $\psi$ be a copy of $J$ in $G$.
	We define $x(\psi), y(\psi)$ to be the number of copies of $H_{i},H_{i-1}$ respectively in $G$ that extend $J$
	(a copy
	$\phi$ of $H_i$ or $H_{i-1}$ {\em extends} $\psi$ if the restriction of
	$\phi$ to $V(J)$ equals $\psi$).
	A copy of $J$ in $G$ is {\em royal} if $x(\psi)\ge f(i)n^2/2$,
	and a copy of $H_i$ in $G$ is {\em noble} if it extends a royal copy of $J$. It follows that at least half of
	the copies of $H_i$ in $G$ are noble.
	
	Let $\psi$ be a royal copy of $J$. We claim that $y(\psi)\ge (\vare/2)x(\psi)$.
	Let $P$ be the set of vertices $v\in V(G)$
	such that mapping $p$ to $v$ extends $\psi$ to a copy of $H_i\setminus \{q\}$, and let $Q$ be
	the set of vertices $v\in V(G)$
	such that mapping $q$ to $v$ extends $\psi$ to a copy of $H_i\setminus \{p\}$. Thus either $P\cap Q=\emptyset$ or $P=Q$.
	
	Suppose first that $P\cap Q=\emptyset$. Since $\psi$ is royal, it follows that $|P|\cdot|Q|\ge x(\psi)\ge f(i)n^2/2$.
	Hence
	there are at least $\vare|P|\cdot|Q|$ nonedges between $P,Q$, from our assumption; but then
	$y(\psi)\ge 
	\vare|P|\cdot|Q|\ge  \vare x(\psi)$ as claimed.
	
	Now suppose that $P=Q$.
	Since $\psi$ is royal,
	there are at least $x(\psi)/2$ edges with both ends in $P$,
	and so $|P|(|P|-1)/2\ge x(\psi)/2\ge f(i)n^2/2$.
	Choose $C\subseteq P$
	with size
	$\lfloor|P|/2\rfloor$, and let $D=P\setminus C$.
	Thus
	$$|C|\cdot |D|\ge |P|(|P|-1)/4\ge x(\psi)/4\ge f(i)n^2/4.$$
	Consequently
	there are at least $\vare|C|\cdot |D|$ nonedges between
	$C,D$, from our assumption; and since each of these gives two ways to extend $J$ to a copy of $H_{i-1}$, and
	$|C|\cdot |D|\ge x(\psi)/4$, it follows that
	$y(\psi)\ge (\vare/2)x(\psi)$.
	
	This proves our claim that $y(\psi)\ge (\vare/2)x(\psi)$, for each royal $\psi$. Summing over all royal $\psi$, we deduce that
	the number of copies of $H_{i-1}$ in $G$ is at least
	$\vare/2$ times the number of noble copies of $H_{i}$ in $G$, and hence at least
	$\vare/4$ times the number of copies of $H_{i}$.  This proves (2).
	
	\bigskip
	
	Since one of $H_0\LL H_m$ equals $H$, and $G$ is $H$-free, it is not the case that $f(i-1)\ge \vare f(i)/4$ for all $i$
	with $1\le i\le m$; choose $i\in \{1\LL m\}$ maximum such that $f(i-1)< \vare f(i)/4$. Consequently,
	$f(j-1)\ge \vare f(j)/4$ for $i+1\le j\le m$, and so by (1),
	$$f(i)\ge (\vare/4)^{m-i}f(m)\ge (\vare/4)^{m-1}t^{-k^2}.$$
	From (2), there are disjoint subsets $W_1,W_2$ of $V(G)$
	with $|W_1|\cdot|W_2|\ge (f(i)/4)n^2$, such that there are fewer than $\vare|W_1|\cdot |W_2|$ nonedges between $W_1,W_2$.
	Since $|W_1|,|W_2|\le n$, it follows that
	$$|W_1|,|W_2|\ge (f(i)/4)n\ge \frac14(\vare/4)^{m-1}t^{-k^2}n\ge 2(\vare/4)^{k(k-1)/2}t^{-k^2}n.$$
	Since there are fewer than $\vare|W_1|\cdot |W_2|$ nonedges between $W_1,W_2$, fewer than half the vertices in $W_1$
	have at least $2\vare|W_2|$ neighbours in $W_2$, and the result follows.~\bbox
	
	Let us deduce \ref{foxsud1} %\ref{foxsudakov} 
	from \ref{foxsudmain}, in the following form (to deduce \ref{foxsudakov} itself, take $\delta=1/2$, say):
	\begin{thm}\label{stronger}
		For every graph $H$, and every $\delta$ with $0<\delta\le 1/2$, let $c=\delta/(6|H|^2)$; then every $H$-free graph $G$
		with $2|G|^{-1}+|G|^{-\delta}\le 1$ contains either a stable set of size $|G|^c$, or a complete bipartite subgraph with
		parts of cardinality at least $|G|^{1-\delta}$
		and at least $|G|^c$ respectively.
	\end{thm}
	\Proof
	Let $H,\delta,c$ be as in the theorem, and let $G$ be $H$-free. Suppose first that $|G|^c\le 2$.
	If $G$ is not complete, the first outcome holds, since $|G|^c\le 2$; and if $G$ is complete, the second holds,
	since $|G|\ge |G|^{1-\delta}+2$. Thus we may assume that $|G|^c>2$.
	Let $t:=\lceil |G|^c\rceil$
	and $\vare :=1/(4t)$.
	
	Since $|G|^c\ge 2$, it follows that
	$$|G|^{\delta-3c|H|^2/2} =|G|^{9c|H|^2/2}\ge 2^{9|H|^2/2} \ge 2^{1+7|H|^2/2}.$$
	Consequently
	$$t^{-|H|^2}\left(\frac{\vare}{4}\right)^{|H|^2/2}=t^{-|H|^2}\left(\frac{1}{16t}\right)^{|H|^2/2}\ge (2|G|^c)^{-3|H|^2/2}4^{-|H|^2}\ge 2|G|^{-\delta}$$
	since $t\le 2|G|^c$.

	We may assume that $G$ has no stable set of size $t$, and so by \ref{foxsudmain}, and since $\delta\le 1/2$ and so $1-\delta\ge c$,
	there are disjoint subsets $W_1,W_2$ of $V(G)$, with
	$$|W_1|,|W_2|\ge t^{-|H|^2}\left(\frac{\vare}{4}\right)^{|H|^2/2}|G|\ge 2|G|^{1-\delta}\ge 2|G|^c\ge t,$$
	such that every vertex in $W_1$ is nonadjacent to at most $2\vare|W_2|$ vertices  in $W_2$.
	Choose $A\subseteq W_1$ with cardinality $t$. Since each vertex in $A$ has at most $2\vare|W_2|$ non-neighbours in $W_1$,
	there are at least
	$|W_2|(1-2\vare t)=|W_2|/2\ge |G|^{1-\delta}$ vertices in $W_2$ adjacent to every vertex in $X$.
	This proves \ref{stronger}.~\bbox
	
	For its application later in this paper, we need a more powerful version of \ref{foxsudmain}, replacing the hypothesis that $G$ is
	$H$-free
	by the weaker hypothesis that $G$ does not contain many copies of $H$, and replacing the stable set outcome
	with a linear sparse set outcome. (Actually, we only need the second strengthening, but the first comes for free anyway.)
	After some preliminaries, its proof is much the
	same as that of \ref{foxsudmain}.
	
	We need a lemma:
	\begin{thm}\label{countcliques}
		Let $h\ge 1$ be an integer and let $0<c<1$. For every graph $G$, either:
		\begin{itemize}
			\item $\ind_{K_h}(G)\ge c^{h(h+1)/2}|G|^h$, or
			\item there exists $S\subseteq V(G)$ with $|S|\ge c^{h(h-1)/2}|G|$ such that $G[S]$ has fewer than $c|S|^2$ edges.
		\end{itemize}
	\end{thm}
	\Proof
	Suppose that $G$ is a graph such that $G[S]$ has at least $c|S|^2$ edges for every $S\subseteq V(G)$ with $|S|\ge c^{h(h-1)/2}|G|$.
	We observe first:
	\\
	\\
	(1) {\em For every set $S\subseteq V(G)$ with $|S|\ge c^{h(h-1)/2}|G|$, there are at least $c|S|$ vertices in $S$ that have degree
		at least $c|S|$ in $G[S]$.}
	\\
	\\
	Let there be $x|S|$ vertices in $S$ that have degree at least $c|S|$. Then the sum of the degrees in $G[S]$ is at most
	$x|S|^2+c|S|^2$, and since $G[S]$ has at least $c|S|^2$ edges, it follows that $x+c\ge 2c$. This proves (1).
	
	\bigskip
	
	Choose $v_1\LL v_h\in V(G)$ independently and uniformly at random. For $1\le i\le h$, let $E_i$ be the event that
	$x_1\LL x_i$ are distinct and pairwise adjacent, and there are at least $c^i|G|$ vertices that are distinct from $x_1\LL x_i$
	and adjacent to all of $x_1\LL x_i$; and let $p_i$ be the probability of $E_i$. We claim that:
	\\
	\\
	(2) {\em $p_i\ge c^{i(i+1)/2}$ for $1\le i\le h$.}
	\\
	\\
	We prove this by induction on $i$. The result holds for $i = 1$, since at least $c|G|$ vertices of $G$ have degree at
	least $c|G|$ by (1). We assume that $i\ge 2$ and the result holds for $i-1$. But $E_i$ is the event that
	\begin{itemize}
		\item $E_{i-1}$
		holds, and
		\item $v_i$ belongs to $X$ where $X$ is the set of vertices that are distinct from $x_1\LL x_{i-1}$
		and adjacent to all of $x_1\LL x_{i-1}$, and
		\item $v_i$ is adjacent to at least $c^i|G|$ members of $X$.
	\end{itemize}
	If $E_{i-1}$ holds, then the set $X$ above has size at least $c^{i-1}|G|\ge c^{h(h-1)/2}|G|$, and so by (1),
	at least $c|X|$ vertices in $X$ are adjacent to at least $c|X|$ vertices in $X$. Consequently
	$$p_i\ge p_{i-1}(c|X|/|G|)\ge p_{i-1}c^i\ge c^{(i-1)i/2}c^i= c^{i(i+1)/2}.$$
	This proves (2).
	
	\bigskip
	
	In particular, the probability that $v_1\LL v_h$ are all distinct and pairwise adjacent is at least $c^{h(h+1)/2}$;
	and so $\ind_{K_h}(G)\ge c^{h(h+1)/2}|G|^h$. This proves \ref{countcliques}.~\bbox
	
	We deduce the following more powerful version of \ref{foxsudmain}:

	\begin{thm}\label{strengthening}
		Let $H$ be a graph with $h\ge 1$ vertices, and let $0<\vare<1/4$. Let $G$ be a graph with $\ind_H(G)<(\vare^{h} |G|)^{h}$.
		Then either
		\begin{itemize}
			\item  there exists $S \subseteq V(G)$ with $|S| \ge \vare^{h(h-1)/2} |G|$ such that $|E(G[S])|<\vare |S|^2$, or
			\item there exist disjoint $W_1,W_2\subseteq V (G)$ with $|W_1|, |W_2| \ge 2(\vare/2)^{h^2}|G|$ such that there are fewer than 
			$\vare|W_1|\cdot |W_2|$ nonedges between $W_1$ and $W_2$.
		\end{itemize}
	\end{thm}
	\Proof Let $n:=|G|$. We may assume there is no $S$ satisfying the first outcome, and so
	$\ind_{K_h}(G)\ge \vare^{h(h+1)/2}n^h$ by \ref{countcliques}. Define $m$ and $H_0,H_1\LL H_m$ as in the proof of \ref{foxsudmain},
	and for $0\le i\le m$, let $f(i)n^h$ be the number of copies of $H_i$ in $G$. Thus we have seen that
	$f(m)\ge \vare^{h(h+1)/2}$.
	As in the proof of \ref{foxsudmain}, we have
	\\
	\\
	(1) {\em For $1\le i\le m$, either $f(i-1)\ge (\vare/4) f(i)$, or there are disjoint subsets $W_1,W_2$ of $V(G)$
		with $|W_1|\cdot|W_2|\ge (f(i)/4)n^2$, such that there are fewer than $\vare|W_1|\cdot |W_2|$ nonedges between $W_1,W_2$.}
	
	\bigskip
	Since one of $H_0\LL H_m$ equals $H$, and $\ind_H(G)<(\vare^{h} |G|)^{h}<(\vare/4)^m$ (because $\vare\le 1/4$),
	it is not the case that $f(i-1)\ge \vare f(i)/4$ for all $i$ with $1\le i\le m$;
	choose $i\in \{1\LL m\}$ maximum such that $f(i-1)< \vare f(i)/4$. Consequently,
	$f(j-1)\ge \vare f(j)/4$ for $i+1\le j\le m$, and so
	$$f(i)\ge (\vare/4)^{m-i}f(m)\ge (\vare/4)^{m-1}\vare^{h(h+1)/2}=\vare^{(h(h-1)/2-1+h(h+1)/2}4^{1-m}=2^{-h^2+h+2}\vare^{h^2-1}\ge 8(\vare/2)^{h^2}.$$
	From (1), there are disjoint subsets $W_1,W_2$ of $V(G)$
	with $|W_1|\cdot|W_2|\ge (f(i)/4)n^2\ge 2(\vare/2)^{h^2}n^2$, such that there are fewer than $\vare|W_1|\cdot |W_2|$ nonedges between $W_1,W_2$.~\bbox
	
	We remark that \ref{strengthening} implies \ref{foxsudmain}, with worse constants: given $t$ and $\vare$ as in \ref{foxsudmain},
	define $\vare'=\min(\vare,1/(2t))$, and apply \ref{strengthening} with $\vare$ replaced by $\vare'$.
	
	%%%%%%%%%%%%%%%%%%%%%%%%%%%%%%%%%%%%%%%%%%%%%%%%%%%%%%%%%%%%%%%%%%%%%%%%%%%%%%%%%%%%%%%%%%%%%%%%%%%%%%%
	
	\section{A sparse sequence of dense pairs}\label{sec:sparseseq}
	%For two graphs $G,H$, and $x>0$, we define \[\mu_H(x,G):=\frac{\ind_H(G)}{(x\abs G)^{\abs H}}.\]
	%For a finite family $\mac F$ of graphs, let $\mu_{\mac F}(x,G):=\max_{H\in\mac F}\mu_H(x,G)$.
	Suppose we are given a sparse $H$-free graph $G$, say with
	density at most $y$.  The goal of this section is to show that either we can
	drop to a much sparser (with density $O(y^2)$) induced subgraph that is still large (size $\poly(y)|G|$), or we can find a `nice' structure: a sequence
	of disjoint sets $A_1,\dots,A_k,B_1,\dots,B_k$ such that pairs $A_i$, $B_i$ are dense
	to each other, while distinct pairs are sparse to each other (and all the
	parameters have suitable polynomial dependence on $y$). We will use these
	structures in later sections to build the induced subgraphs (Swiss army
	graphs) that we are looking for.
	
	If $A,B\subseteq V(G)$ are disjoint, and $0\le x\le 1$, we say that $B$ is {\em $x$-sparse} to $A$ if every vertex in $B$ has at most $x|A|$
	neighbours in $A$; and $B$ is {\em $x$-dense} to $A$ if every vertex in $B$ has at least $x|A|$ neighbours in $A$.
	\begin{thm}
		\label{lem:sparsedense}
		Let $H$ be a graph, let
		$y\in(0,2^{-99})$, and let $G$ be an $H$-free $y$-sparse graph.
		Suppose that 
		there is no subset $S\subseteq V(G)$ with $|S| \ge y^{|H|^2} |G|$ such that $G[S]$ is $y^2$-sparse.
		Then there exist disjoint $A_1,A_2,\ldots,A_k,B_1,B_2,\ldots,B_k\subseteq V(G)$, where $k=\ceil{y^{-1/4}}$, such that:
		\begin{itemize}
			\item $\abs{A_i},\abs{B_i}= \ceil{ (y^2/16)^{|H|^2}|G|}$, and 
			there are at most $(y^2/2)\abs{A_i}\abs{B_i}$ nonedges between $A_i,B_i$, for $1\le i\le k$; and
			\item 
			each of $A_i,B_i$ is $y^{1/6}$-sparse to each of $A_j,B_j$, for all distinct $i,j\in \{1\LL k\}$.
		\end{itemize}
	\end{thm}
	\Proof Let $h:=|H|$. Since there is no $S$ as in the theorem, it follows that $E(G)\ne \emptyset$ and so  $y|G|\ge 1$.
	Let $s\ge0$ be maximum such that there are disjoint $X_1,\ldots,X_s,Y_1,\ldots,Y_s\subseteq V(G)$ satisfying:
	\begin{itemize}
		\item for all $i\in[s]$, $\abs{X_i},\abs{Y_i}=\ceil{ 2(y^2/16)^{h^2}|G|}$, and there are at most $(y^2/8)\abs{X_i}\abs{Y_i}$ nonedges between $X_i,Y_i$; and
		
		\item for $1\le i<j\le s$, each of $X_j,Y_j$ is $y^{1/2}$-sparse to each of $X_i,Y_i$.
	\end{itemize}
	We claim that:
	\\
	\\
	(1) {\em $s\ge y^{-1/4}$.}
	\\
	\\
	Suppose not. For $1\le i\le s$,
	$$|X_i|=|Y_i|=\ceil{ 2(y/16)^{h^2}|G|}\le \max(1, 4(y/16)^{h^2}|G|)\le y|G|$$ 
	(since $y|G|\ge 1$).
	Let $A:=\bigcup_{i\in [s]}(X_i\cup Y_i)$.
	For each $i\in[s]$, let $D_i$ be the set of vertices in $V(G)\setminus A$ that have either at least 
	$y^{1/2}\abs{X_i}$ neighbours in $X_i$ or at least $y^{1/2}\abs{Y_i}$ neighbours in $Y_i$;
	then $\abs{D_i}\le 2y^{1/2}\abs G$, since there are at most $y\abs{X_i}|G|$ edges between $X_i$ and $V(G)\setminus A$,
	and the same for $Y_i$. 
	Let $G':=G\setminus(A\cup\bigcup_{i\in[s]}D_i)$; then since $s<y^{-1/4}$, we have
	\[\abs{G'}\ge\abs G-s(2y\abs G+2y^{1/2}\abs G)
	\ge \abs G-2y^{-1/4}(y\abs G+y^{1/2}\abs G)
	\ge (1-4y^{1/4})\abs G\ge |G|/2.\]
	By \ref{strengthening} applied to $G'$, taking $\vare=y^2/8$,
	either
	\begin{itemize}
		\item  there exists $S \subseteq V(G')$ with $|S| \ge (y^2/8)^{h(h-1)/2} |G'|\ge y^{h^2}|G|$ such that $|E(G[S])|<(y^2/8)|S|^2$, or
		\item there exist disjoint $W_1,W_2\subseteq V (G')$ with $|W_1|, |W_2| \ge 2(y^2/16)^{h^2}|G'|$ such that there are fewer than
		$(y^2/8)|W_1|\cdot |W_2|$ nonedges between $W_1$ and $W_2$. 
	\end{itemize}
	Suppose that the first holds, and let $S$ be the corresponding set. Since $|E(G[S])|<(y^2/8)|S|^2$, at most $|S|/2$
	vertices in $S$ have degree in $G[S]$ at least $y^2|S|/2$; and so there is a subset $S'\subseteq S$  that is $y^2$-sparse, with 
	$|S'|\ge |S|/2\ge (y^2/8)^{h(h-1)/2} |G'|/2\ge y^{|H|^2} |G|$, contrary to the hypothesis.
	
	If the second bullet holds, then by averaging, there are subsets $X_{s+1}\subseteq W_1$ and
	$Y_{s+1}\subseteq W_2$, both of size $\ceil{ 2(y^2/16)^{h^2}|G|}$, such that there are fewer than $(y^2/8)|X_{s+1}|\cdot |Y_{s+1}|$ nonedges between
	$X_{s+1}$ and $Y_{s+1}$,
	contrary to the maximality of $s$. This proves (1).
	
	\bigskip
	
	Let $k=\lceil y^{-1/4}\rceil$.
	For $1\le i\le k$, and for each $j$ with $i<j\le k$, there are at most $y^{1/2}|X_i|\cdot |X_j|$ edges between $X_i,X_j$ (since $X_j$
	is $y^{1/2}$-sparse to $X_i$), and so at most $2y^{1/3}\abs{X_i}$ vertices in $X_i$ have at least $\frac12y^{1/6}\abs{X_j}$ 
	neighbours in $X_j$. Similarly there are at most $2y^{1/3}\abs{X_i}$ vertices in $X_i$ that have at least 
	$\frac12y^{1/6}\abs{Y_j}$ neighbours in $Y_j$. Let $A_i'\subseteq X_i$ be the set of vertices $v\in X_i$ such that for each 
	$j\in \{i+1\LL k\}$, $v$ has at most $\frac12y^{1/6}\abs{X_j}$ neighbours in $X_j$ and at most 
	$\frac12y^{1/6}\abs{Y_j}$ neighbours in $Y_j$. It follows that $|A_i'|\ge (1-2y^{1/3}y^{-1/4})|X_i|\ge |X_i|/2$; and so there
	exists $A_i\subseteq A_i'$ with $|A_i|=\ceil{|X_i|/2}= \ceil{ (y^2/16)^{h^2}|G|}$.
	Define $B_i\subseteq Y_i$ similarly: that is, $B_i$ is a  set of vertices $v\in Y_i$ such that for each $j\in \{i+1\LL k\}$, $v$ has at most $\frac12y^{1/6}\abs{X_j}$ neighbours in $X_j$ and at most
	$\frac12y^{1/6}\abs{Y_j}$ neighbours in $Y_j$, and $|B_i|= \ceil{ (y^2/16)^{h^2}|G|}$. 
	
	Consequently there are at most $(y^2/8)\abs{X_i}\abs{Y_i}\le (y^2/2)\abs{A_i}\abs{B_i}$ nonedges between $A_i,B_i$.
	Now, for all $i,j\in \{1\LL k\}$ with $i<j$,
	each of $A_i,B_i$ is $y^{1/6}$-sparse to each of $A_j,B_j$, since they are $\frac12y^{1/6}$ sparse to $X_j$ and to $Y_j$;
	and each of $A_j,B_j$ is $y^{1/6}$-sparse to each of $A_i,B_i$ since each of $A_j,B_j$ is $y^{1/2}$-sparse to each of 
	$X_i,Y_i$ and $2y^{1/2}\le y^{1/6}$.
	This proves \ref{lem:sparsedense}.~\bbox
	
	%%%%%%%%%%%%%%%%%%%%%%%%%%%%%%%%%%%%%%%%%%%%%%%%%%%%%%%%%%%%%%%%%%%%%%%%%%%%%%%%%%%%%%%%%%%%%%%%%%%%%%
	\section{Blockades, and growing a hand}\label{sec:blockades}

	A {\em blockade} in $G$ is a sequence $\mathcal{B}=(B_1\LL B_k)$ of
	pairwise disjoint
	subsets of $V(G)$, and we call $B_1\LL B_k$ its {\em blocks}. (In some earlier papers, the blocks of a blockade must be nonempty,
	but here it is convenient to allow empty blocks.) We define $V(\mac B)=B_1\cupcup B_k$
	The {\em length} of the blockade $\mathcal{B}=(B_1\LL B_k)$ is $k$, and its {\em width} is the minimum of the cardinalities
	of its blocks. If its length is at least $\ell$ and width at least $w$ we call it an {\em $(\ell,w)$-blockade}.
	For $\vare>0$, the blockade
	$\mathcal{B}=(B_1\LL B_k)$ is
	{\em $\vare$-sparse} if $B_{i+1}\cupcup B_k$ is $\vare$-sparse to $B_i$ for all $i$ with $1\le i\le k$; and similarly
	$\mathcal{B}$ is {\em $(1-\vare)$-dense} if  $B_{i+1}\cupcup B_k$ is $(1-\vare)$-dense to $B_i$ for all $i$ with $1\le i\le k$.

	Let $\alpha(G), \omega(G)$ be the cardinalities of the largest stable sets and cliques 
	of $G$ respectively; and for $c>0$, let us say a graph $G$ is {\em $c$-critical} if $\alpha(G)\omega(G)<|G|^c$, 
	and 
	$\alpha(G')\omega(G')\ge |G'|^c$ for every induced subgraph
	$G'$ of $G$ with $|G'|<|G|$. In order to prove \ref{swiss}, it suffices to show that for every two Swiss Army graphs $H,J$,
	if $c>0$ is sufficiently small, then no $\{H,\overline{J}\}$-free graph is $c$-critical. 
	If $X,Y\subseteq V(G)$ are disjoint, we say $X$ is {\em anticomplete} to $Y$ if there are no edges of $G$ between $X$ and $Y$.
	
	As a first step, we will prove that if $c$ is sufficiently small, then every $c$-critical graph contains a piece of machinery
	that will provide the 2-subdivision paths of a Swiss Army graph. Thus, let $\mathcal{B}=(B_1\LL B_k)$ be a blockade in $G$. 
	A {\em $B_i$-finger} is an induced path in $G$ of length two, with three vertices $a_i\DD b_i\DD c_i$ in order, such that
	\begin{itemize}
		\item $a_i,b_i,c_i\notin V(\mac B)$;
		\item $c_i$ is complete to $B_i$ and is anticomplete to $V(\mac B)\setminus B_i$; and
		\item $a_i,b_i$ are both anticomplete to $V(\mac B)$.
	\end{itemize}
	Let $a_i\DD b_i\DD c_i$ be a $B_i$-finger, for $1\le i\le k$. We call the union of these fingers a {\em hand} for $\mathcal{B}$
	if 
	\begin{itemize}
		\item $a_1,a_2\LL a_k$ are all the same vertex;
		\item $b_1,c_1,b_2,c_2\LL b_k,c_k$ are all distinct;
		\item the union of these fingers is induced; that is, for $1\le i<j\le k$, $\{b_i,c_i\}$ is anticomplete to $\{b_j,c_j\}$.
	\end{itemize}
	We call $a_1$ the {\em palm} of the hand.
	
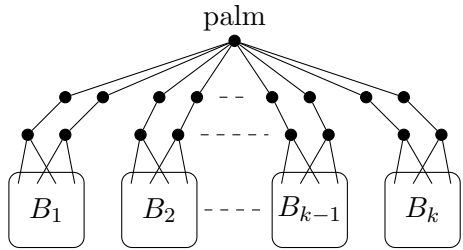
\begin{figure}[h!]
		\centering
		
		\begin{tikzpicture}[scale=1/2,auto=left]
			
			\tikzstyle{every node}=[inner sep=1.5pt, fill=black,circle,draw]
			
			\draw[rounded corners] (-6,-1) rectangle (-4,1);
			\draw[rounded corners] (-3,-1) rectangle (-1,1);
			\draw[rounded corners] (1,-1) rectangle (3,1);
			\draw[rounded corners] (4,-1) rectangle (6,1);
			
			\node (c1) at (-5.5,2) {};
			\node (c2) at (-4.5,2) {};
			\node (c3) at (-2.5,2) {};
			\node (c4) at (-1.5,2) {};
			\node (c5) at (1.5, 2) {};
			\node (c6) at (2.5,2) {};
			\node (c7) at (4.5,2) {};
			\node (c8) at (5.5,2) {};
			
			\draw (c1) to (-5.75, .7);
			\draw (c2) to (-5.25, .7);
			\draw (c1) to (-4.75, .7);
			\draw (c2) to (-4.25, .7);
			
			\draw (c3) to (-2.75, .7);
			\draw (c4) to (-2.25, .7);
			\draw (c3) to (-1.75, .7);
			\draw (c4) to (-1.25, .7);
			
			\draw (c5) to (1.25, .7);
			\draw (c6) to (1.75, .7);
			\draw (c5) to (2.25, .7);
			\draw (c6) to (2.75, .7);
			
			\draw (c7) to (4.25, .7);
			\draw (c8) to (4.75, .7);
			\draw (c7) to (5.25, .7);
			\draw (c8) to (5.75, .7);
			
			\node (b1) at (-4.5,3) {};
			\node (b2) at (-3.5,3) {};
			\node (b3) at (-2,3) {};
			\node (b4) at (-1,3) {};
			\node (b5) at (1, 3) {};
			\node (b6) at (2,3) {};
			\node (b7) at (3.5,3) {};
			\node (b8) at (4.5,3) {};
			
			\node (a) at (0,4.5) {};
			
			\foreach \from/\to in {a/b1,a/b2,a/b3,a/b4,a/b5,a/b6,a/b7,a/b8,b1/c1,b2/c2,b3/c3,b4/c4,b5/c5,b6/c6,b7/c7,b8/c8}
			\draw [-] (\from) -- (\to);
			
			\draw[dashed] (-.8,0) -- (.8,0);
			\draw[dashed] (-.9,2) -- (.9,2);
			\draw[dashed] (-.4,3) -- (.4,3);
			
			\tikzstyle{every node}=[]
			\draw (a) node [above]           {palm};
			\node at (-5,0) {$B_1$};
			\node at (-2,0) {$B_2$};
			\node at (2,0) {$B_{k-1}$};
			\node at (5,0) {$B_k$};
			
		\end{tikzpicture}
		
		\caption{A 2-thickened hand for $(B_1\LL B_k)$.} \label{fig:thickhand}
	\end{figure}

	Now let $H_1\LL H_s$ be hands for $\mac B$, all with the same palm but otherwise pairwise vertex-disjoint and anticomplete, 
	and take their union. We call the result an {\em $s$-thickened hand}, with {\em palm} the common palm of $H_1\LL H_s$.
	See figure \ref{fig:thickhand}.
	
	Finally, let $T_1\LL T_t$ be $s$-thickened hands for $\mathcal{B}$, pairwise disjoint and anticomplete. Take their union,
	and for $1\le i<j\le t$ add $s$ paths of length three joining the palms of $T_i,T_j$ (that is, add $s$ edges each joining the two given palms, and then subdivide each of them twice).
	The graph we produce is called
	a {\em $(s,t)$-handset} for $\mathcal B$. 
	
	We need the following easy lemma:
	\begin{thm}\label{getstable}
		Let $F$ be a graph and let $t,n\ge 0$ and $m\ge 1$ be integers, with $|F|\ge m^tn$. Then either:
		\begin{itemize}
			\item $F$ has a stable set of size $m$; or
			\item there are disjoint subsets $X,Y$ of $V(F)$, such that $X$ is a clique and $|X|=t$,
			and $X$ is complete to $Y$ and $|Y|\ge n$.
		\end{itemize}
	\end{thm}
	\Proof
	We proceed by induction on $t$. If $t=0$ the result is true (take $X=\emptyset$ and $Y=V(F)$), so we assume that $t>0$ and the result 
	holds for $t-1$. We may assume that $F$ has no stable set of size $m$, and so its chromatic number is more than $|F|/m$,
	and therefore it has a vertex with more than $|F|/m-1\ge m^{t-1}n-1$ neighbours, and hence with at least $m^{t-1}n$ neighbours.
	The result follows from the inductive hypothesis applied to the subgraph induced on this set of neighbours. This proves \ref{getstable}.~\bbox
	
	A blockade $(B_1\LL B_k)$ is {\em equicardinal} if all its blocks have the same size; and {\em symmetrically $x$-sparse}
	if $B_i$ is $x$-sparse to $B_j$ for all distinct $i,j\in \{1\LL k\}$.
	We will prove the following.
	
	\begin{thm}\label{gethandset}
		Let $s,t\ge 1$ be integers, let $\rho>0$, and let $H$ be a graph. Then there exist $d,\delta,\eta>0$ with the following property. Let $0\le c\le \delta$
		,
		and let $G$ be a $c$-critical $H$-free graph. Let $Z\subseteq V(G)$ with $|Z|\ge y^\rho |G|$ such that $G[Z]$ is $y$-sparse for some $y$ with $0<y<\eta$.
		Then either:
		\begin{itemize}
			\item there is a subset $S\subseteq Z$ with $|S| \ge y^{|H|^2} |Z|$ that is $y^2$-sparse; or
			\item there is a symmetrically $2y^{1/6}$-sparse equicardinal blockade $\mathcal{B}$ in $G$ of length at least $y^{-1/64}$ and width at least 
			$y^d|Z|/2$, and there is an $(s,t)$-handset for $\mac B$.
		\end{itemize}
	\end{thm}
	\Proof
	
	%$k=\ceil{y^{-1/L}}$
	%sk^2y^{1/6} \le 1/2
	%$m (y^d/2)^c> 2$
	%$k\ge  m^tn$
	%$2ys(t+n) + (2st(t+n)+t) y^{1/6}<1/2$
	Let $h:=|H|$.
	Choose $\eta>0$ such that
	$$\eta\le \min\left( 2^{-99}, \left(16st^2\right)^{-6}, (8s)^{-24}, 2^{-96t}\right)$$
	and 
	$$\eta s\left(t+2\eta^{-1/64}\right) + \left(2st\left(t+2\eta^{-1/64}\right)+t\right) \eta^{1/6}<1/2.$$
	Choose $d>1$ such that $(y^2/16)^{h^2}\ge y^d$ for all $y\le \eta$, and let $\delta=1/(64(d+\rho)t)$.
	We claim that $d, \delta, \eta$ satisfy the theorem.

	Let $0<c\le \delta$, let $G$ be a $c$-critical $H$-free graph, and let $Z\subseteq V(G)$ with $|Z|\ge y^\rho$, such that 
	$G[Z]$ is $y$-sparse,
	where $0<y<\eta$.
	We may assume that the first outcome of the theorem, applied to $G[Z]$, is false. So
	by \ref{lem:sparsedense}, 
	there exist disjoint 
	$$A_1,A_2,\ldots,A_{k'},B_1,B_2,\ldots,B_{k'}\subseteq Z,$$ 
	where $k'=\ceil{y^{-1/4}}$, such that:
	\begin{itemize}
		\item $\abs{A_i},\abs{B_i}= \ceil{ (y^2/16)^{h^2}|Z|}$, and
		there are at most $(y^2/2)\abs{A_i}\cdot \abs{B_i}$ nonedges between $A_i,B_i$, for $1\le i\le k'$; and
		\item
		each of $A_i,B_i$ is $y^{1/6}$-sparse to each of $A_j,B_j$, for all distinct $i,j\in \{1\LL k'\}$.
	\end{itemize}
	Let $k:= \ceil{y^{-1/16}}$; we will only need the sets $A_i,B_i$ for $1\le i\le k$.
	
	Let $1\le i\le k$. Since there are at most $(y^2/2)\abs{A_i}\cdot \abs{B_i}$ nonedges between $A_i,B_i$, there exists $X\subseteq B_i$
	with $|X|\ge |B_i|/2$ that is $(1-y^2)$-dense to $A_i$. Since $G$ is $c$-critical, the graph $G[X]$ satisfies
	$\alpha(G[X])\omega(G[X])\ge |X|^c$; and since $\omega(G[X])\le \omega(G)$, we deduce that there is a stable subset $C_i\subseteq B_i$
	of size at least $(|B_i|/2)^c/\omega(G)$ that is $(1-y^2)$-dense to $A_i$.

	Let $C:=C_1\cupcup C_k$. An {\em induced matching} of $G[C]$ means a set of edges of $G[C]$, pairwise vertex-disjoint and anticomplete.
	Choose an induced matching $M$ of $G[C]$, maximal with the property that for all distinct $i,j\in \{1\LL k\}$, at most $s$ members of $M$
	have an end in $C_i$ and an end in $C_j$.
	It follows that $|M|\le sk^2/2$.
	Let $F$ be the graph with vertex set $\{1\LL k\}$ in which $i,j$ are adjacent if there are $s$ members of $M$ with an end in $C_i$
	and an end in $C_j$. Let $m:=\ceil{ (y^2/16)^{-ch^2}y^{-c\rho}2^{c}}$, and $n:=\ceil{y^{-1/64}}$.
	Thus 
	$$m\le 2(y^2/16)^{-ch^2}y^{-c\rho}2^{c}\le y^{-cd-c\rho}2^{c+1}.$$
	
	%$m=\ceil{ (y^2/16)^{-ch^2}y^{-c\gamma}2^{c}}$
	By \ref{getstable}, either:
	\begin{itemize}
		\item $k< m^tn$; or
		\item $F$ has a stable set of size $m$; or
		\item there are disjoint subsets $X,Y$ of $V(F)$, such that $X$ is a clique of $F$ and $|X|=t$,
		and $X$ is complete to $Y$ in $F$, and $|Y|\ge n$.
	\end{itemize}
	
	Suppose the first bullet holds. Since $m<y^{-cd-c\rho}2^{c+1}$, and $n\le 2y^{-1/64}$, it follows that
	$$y^{-1/16}\le k< 2^{t(c+1)}y^{-c(d+\rho)t}\left(2y^{-1/64}\right).$$
	Since $c(d+\rho)t\le 1/64$, we deduce that $y^{-1/16} <2^{t(c+1)+1}y^{-1/32}$, that is, $y^{-1/32}<2^{t(c+1)+1}<2^{3t}$, a contradiction
	since $y<2^{-96t}$.
	%y<2^{-128t}$

	Suppose the second bullet holds; then we may assume that $I$ is a stable set of $F$ and $|I|=m$. 
	Let $P$ be the set of vertices in 
	$\bigcup_{i\in I}C_i$ that either belong to a member of $M$ or are adjacent to a member of $M$. Thus for each $i\in I$, 
	$$|P\cap C_i|\le 2y^{1/6}|M|\cdot |C_i|\le sk^2y^{1/6}|C_i|\le |C_i|/2,$$ 
	since $k\le 2y^{-1/16}$, and $y\le (8s)^{-24}$. But from the maximality of $M$, the union of the sets 
	%sk^2y^{1/6} \le 1/2
	$C_i\setminus P\;(i\in I)$ is stable in $G$, and so has cardinality at most $|G|^c/\omega(G)$. Consequently
	$$\sum_{i\in I}|C_i|/2< |G|^c/\omega(G).$$
	But $|C_i|\ge (|B_i|/2)^c/\omega(G)$ for each $i\in I$, and so
	$$\sum_{i\in I}(|B_i|/2)^c/\omega(G)< |G|^c/\omega(G).$$
	Each $|B_i|\ge (y^2/16)^{h^2} |Z|$, and so 
	$m (y^2/16)^{ch^2}(|Z|/2)^c< |G|^c$. Since $|Z|\ge y^\rho |G|$, it follows that
	$m (y^2/16)^{ch^2}y^{c\rho}2^{-c}< 1$, a contradiction, from the choice of $m$.
	%$m (y^d/2)^c> 2$
	
	Suppose the third bullet holds, and let $X,Y$ be the corresponding subsets of $V(F)$. For each edge $ij$ of $F$ with one end in $X$
	and the other end in $X\cup Y$, there  are $s$ edges in $M$ between $C_i$ and $C_j$; let $M'$ be the set of all these edges.
	Thus $|M'|\le st(t+n)$. Let $N$ be the set of ends of all the edges in $M'$. 
	For each $i\in X$, choose $a_i\in A_i$ as follows: $a_i$ is adjacent to each vertex in $C_i\cap N$,
	and nonadjacent to every vertex in $N\setminus C_i$, and the vertices $a_i\;(i\in I)$ are pairwise nonadjacent. 
	To see that this is possible,
	observe that, since $C_i$ is $(1-y)$-dense to $A_i$, and $|N\cap C_i|\le s(t+n)$, there are at most $y|A_i|s(t+n)$
	vertices in $A_i$ that have a non-neighbour in $C_i\cap N$; and since $C\setminus C_i$ is $y^{1/6}$-sparse to $A_i$,
	there are at most $(2st(t+n)+t) y^{1/6}|A_i|$ vertices in $A_i$ that have a neighbour in $N\setminus C_i$ or are adjacent to some 
	already-selected $a_j$. Since
	$$y|A_i|s(t+n) + (2st(t+n)+t) y^{1/6}|A_i|<|A_i|/2$$ 
	from the definition of $\eta$, such a choice of $a_i$ is possible.
	%$2ys(t+n) + (2st(t+n)+t) y^{1/6}<1/2$
	
	For each $j\in J$, let $D_j'$ be the set of vertices in $A_j$ that are adjacent to every vertex in $N\cap C_j$, and nonadjacent
	to every vertex in $N\setminus C_j$, and nonadjacent to the vertices $a_i\;(i\in X)$. By the same argument, $|D_j'|\ge |A_j|/2$
	for each $j\in Y$. Choose $D_j\subseteq D_J'$ of size $\ceil{|A_j|/2}$. But then $(D_j:j\in Y)$ is an equicardinal  
	$(n,y^d\abs G /2)$-blockade, that is symmetrically $2y^{1/6}$-sparse;
	and there is an $(s,t)$-handset for it. This proves \ref{gethandset}.~\bbox
	
	%%%%%%%%%%%%%%%%%%%%%%%%%%%%%%%%%%%%%%%%%%%%%%%%%%%%%%%%%%%%%%%%%%%%%%%%%%%%%%%%%%%%%%%%
	\section{From a sparse blockade to an anticomplete blockade}\label{sec:sparsetoanti}
	
	A blockade is {\em anticomplete} if every two of its blocks are anticomplete.
	If $\mathcal{B}=(B_1\LL B_k)$ is a blockade in $G$, we say an induced subgraph $H$ of $G$ is {\em $\mathcal{B}$-rainbow} if $V(H)\subseteq V(\mathcal{B})$ and
	$|B_i \cap V(H)|\le 1$ for $1\le i\le k$.
	In this section we prove the following (and then we apply it to complete to the proof of \ref{swiss}):
	\begin{thm}\label{getanti}
		Let $F$ be a forest, and let $\alpha, \beta>0$. Then there exist $\alpha',\beta'>0$ with the following property.
		Suppose that $0<y\le 1$, and $G$ is a graph, and $\mac B$ is a $(y^{-\alpha},y^\beta|G|)$-blockade in $G$ that is equicardinal and symmetrically
		$y^\alpha$-sparse. If there is no $\mac B$-rainbow copy of $F$,
		then $G$ admits an anticomplete $(y^{-\alpha'},y^{\beta'}|G|)$-blockade.
	\end{thm}
	We need the following theorem of~\cite{pure1}:
	
	\begin{thm}\label{pure1}
		For every forest $F$, there is an integer $d>0$ with the following property. Let $G$ be a graph
		with a blockade $\mathcal{B}$ of length at least $d$, and let $w$ be the width of $\mathcal{B}$. If
		every vertex of $G$ has degree less than $w/d$, and
		there is no anticomplete pair $A,B\subseteq V(G)$ with $|A|,|B|\ge w/d$,
		then
		there is a $\mathcal{B}$-rainbow copy of $F$ in $G$.
	\end{thm}
	
	This implies:
	
	\begin{thm}\label{firststep}
		For every forest $F$, let $d$ be as in \ref{pure1}. Let $G$ be a graph
		with a blockade $\mathcal{B}$ of length at least $3d^2$, and let $w$ be the width of $\mathcal{B}$. If
		$\mathcal{B}$ is equicardinal and symmetrically $1/d^2$-sparse, and
		there is no $\mathcal{B}$-rainbow copy of $F$, then there is an anticomplete pair $A,B\subseteq V(G)$ with $|A|,|B|\ge w$.
	\end{thm}
	\Proof
	Let $G$ be a graph with a symmetrically $1/d^2$-sparse equicardinal blockade $\mathcal{B}=(B_1\LL B_D)$ of width $w$, where $D=3d^2$. 
	Let $G'$ be the subgraph with vertex set $V(\mac B)$ and edge set the edges of $G$ that have ends that belong to different blocks 
	of $\mac B$. Thus $G'$ has maximum degree at most $(D-1)w/d^2\le 3w$.
	Partition $\{1\LL D\}$
	into $d$ sets of cardinality $3d$, say $I_1\LL I_{d}$. Let $B_h'=\bigcup_{i\in I_h}B_i$ for $1\le i\le d$; then
	$\mathcal{B}'=(B_1'\LL B_d')$ is a $(d,3wd)$-blockade. We may assume there is no $\mac{B}'$-rainbow copy of $F$;
	and so by \ref{pure1} applied to $\mathcal{B}'$, there is an anticomplete pair $(A,B)$ with $|A|,|B|\ge 3w$.
	Choose $i\in \{1\LL k\}$ minimum such that one of 
	$$A\cap (B_1\cupcup B_i), B\cap (B_1\cupcup B_i)$$
	has cardinality at least $w$, and we may assume that $|A\cap (B_1\cupcup B_i)|\ge w$. From the minimality of $i$,
	$|B\cap (B_1\cupcup B_{i-1})|<w$, and since $|B_i|=w$ and $|B|\ge 3w$, it follows that $|B\cap (B_{i+1}\cupcup B_k)|\ge w$.
	But then $A\cap (B_1\cupcup B_i),B\cap (B_{i+1}\cupcup B_k)$ is a pair of subsets of $V(G)$ that are anticomplete
	in $G$ (not just in $G'$), and both have size at least $w$. This proves \ref{firststep}.~\bbox

	\begin{thm}\label{blockparty}
		Let $F$ be a forest, and let $d$ be as in \ref{pure1}; then
		for every integer $s\ge 1$ and every graph $G$, the following holds. Let $D=2 (2d^2)^{s}$, and let $\mathcal{B}$ be a
		%D_1\ge 3d^2
		blockade in $G$ of length $D$, that is equicardinal and symmetrically $2/(2d^2)^s$-sparse, such that there is no 
		$\mathcal{B}$-rainbow copy of $F$. Then
		$G$ admits an anticomplete $(2^s, w/(2d^2)^{s-1})$-blockade,
		where $w$ is the width of $\mathcal{B}$.
	\end{thm}
	\Proof This is true if $s=1$, from \ref{firststep}, and so we assume it is true for some $s-1\ge 1$ and prove it for $s$.
	Let $D=2(2d^2)^{s}$, and let $G$ be a graph with an equicardinal, symmetrically $2/(2d^2)^{s}$-sparse blockade 
	$\mathcal{B}=(B_1\LL B_D)$ 
	of width $w$, and there is no $\mathcal{B}$-rainbow copy of $F$. Partition $\{1\LL D\}$
	into $d$ sets of cardinality $D/d$, say $I_1\LL I_d$. Let $B_h'=\bigcup_{i\in I_h}B_i$ for $1\le i\le d$; then
	$\mathcal{B}'=(B_1'\LL B_d')$ is an equicardinal, symmetrically $2/(2d^2)^{s}$-sparse $(d,wD/d)$-blockade. 
	Let $G'=G[B_1\cupcup B_D]$.
	It follows that there is no $\mathcal{B}'$-rainbow copy of
	$F$; so from \ref{pure1}, there is an anticomplete pair $(A,B)$ of $G'$ with $|A|,|B|\ge wD/d^2$.
	
	Let $D'=2(2d^2)^{s-1}$, and 
	let $I$ be the set of all $i\in \{1\LL D\}$ such that $|A\cap B_i|\ge w/(2d^2)$. Then
	$|I|w+Dw/(2d^2)\ge |A|\ge wD/d^2$, and so $|I|\ge D/(2d^2)=D'$. 
	For each $i\in I$ choose $C_i\subseteq A\cap B_i$ of cardinality
	$\ceil{w'}$, and let $\mathcal{C}$ be the blockade $(C_i:i\in I)$.
	Then $\mathcal{C}$ is equicardinal, of width at least $w/(2d^2)$, and it is symmetrically $2/(2d^2)^{s-1}$-sparse, and 
	there is no $\mathcal{C}$-rainbow copy of
	$F$. Thus the inductive hypothesis, applied to $\mathcal{C}$,
	implies that $G[A]$ admits an
	anticomplete blockade, of width at least $(w/(2d^2))/(2d^2)^{s-2}=w/(2d^2)^{s-1}$ and length $2^{s-1}$; and similarly so does $G[B]$.
	But then combining these gives an anticomplete $(2^s,w/(2d^2)^{s-1})$-blockade in $G$.
	This proves \ref{blockparty}.~\bbox

	We apply this to prove the following strengthened form of \ref{getanti} (in which logarithms are to base two):
	\begin{thm}\label{getanti2}
		Let $F$ be a forest, and let $\alpha, \beta>0$. Let $d$ be as in \ref{pure1}, with $d\ge 8$. Define $\alpha'=\alpha/(5\log d)$ and $\beta'=\alpha+\beta$.
		Suppose that $0<y\le 1$, and there is a $(y^{-\alpha},y^\beta |G|)$-blockade $\mac B$ in a graph $G$ that is equicardinal and symmetrically 
		$y^\alpha$-sparse. If there is no $\mac B$-rainbow copy of $F$,
		then $G$ admits an anticomplete $(y^{-\alpha'}, y^{\beta'} |G|)$-blockade.
	\end{thm}
	\Proof
	Choose an integer $s$ maximal such that $y^{-\alpha}\ge 2(2d^2)^s$. It follows that $s\ge 1$, and 
	$$y^{-\alpha}\le 2(2d^2)^{s+1}\le d^{5s}$$ 
	(since $d\ge 8$ and therefore $2^{s+2}\le d^s$). Thus $\mac B$ has length at least $2(2d^2)^s$ and is equicardinal and symmetrically 
	$2/(2d^2)^2$-sparse. By \ref{blockparty}, $G$ admits an anticomplete $(2^s,w/(2d^2)^{s-1})$-blockade,
	where $w$ is the width of $\mathcal{B}$. But 
	$$2^s= d^{5s/(5\log d)}\ge y^{-\alpha/(5\log d)}=y^{-\alpha'},$$
	and 
	$$w/(2d^2)^{s-1}\ge w/(2(2d^2)^s)\ge y^\alpha w\ge y^{\alpha'} |G|.$$
	This proves \ref{getanti2}.~\bbox
	
	\begin{thm}\label{useanticom} 
		Let $F$ be a forest, and let $\alpha, \beta, \gamma>0$. Then there exists $\delta'>0$ such that for all $c$ with $0<c\le \delta'$, if $G$ is $c$-critical, and 
		$\mathcal{B}$ 
		is a 
		$(y^{-\alpha},y^\beta |G|)$-blockade in $G$ that is equicardinal and symmetrically 
		$y^\gamma$-sparse, then there is a $\mathcal{B}$-rainbow copy of $F$.
	\end{thm}
	\Proof By reducing $\alpha$ or $\gamma$, we may assume that $\alpha=\gamma$ without loss of generality.
	Choose $\alpha',\beta'$ as in \ref{getanti}. 
	Choose $\delta'>0$ such that $\delta'<\alpha'/\beta'$.
	We claim that $\delta'$ satisfies the theorem. Let $0<c\le \delta'$, and let  $\mathcal{B}$
	be a
	$(y^{-\alpha},y^\beta |G|)$-blockade in a $c$-critical graph $G$, that is equicardinal and symmetrically
	$y^\alpha$-sparse, where $0<y\le 1$. Suppose there is no $\mathcal{B}$-rainbow copy of $H$. By \ref{getanti}, 
	$G$ admits an anticomplete $(y^{-\alpha'},y^{\beta'} |G|)$-blockade $\mathcal{A}$.
	Let $\mac A = (A_1\LL A_k)$ say. Since $G$ is $c$-critical, for $1\le i\le k$ there exists a stable subset $C_i\subseteq A_i$ with 
	$|C_i|\omega(G[A_i])\ge |A_i|^c$, and hence with $|C_i|\ge y^{\beta' c}|G|^c/\omega(G)$.
	The the union $C_1\cupcup C_k$ is stable, and hence has cardinality less than $|G|^c/\omega(G)$; and so, since $k\ge y^{-\alpha'}$,
	it follows that 
	$$y^{-\alpha'}y^{\beta' c}|G|^c/\omega(G)\le |G|^c/\omega(G),$$
	that is, 
	$y^{\beta' c-\alpha'}\le 1 ,$
	a contradiction, since $\beta' c-\alpha'<0$. This proves \ref{useanticom}.~\bbox

	We use this to complete the proof of \ref{swiss}, which we restate:
	\begin{thm}\label{swiss2}
		If $H,J$ are Swiss Army graphs, then $\{H,\overline{J}\}$ has the Erd\H{o}s-Hajnal property.
	\end{thm}
	\Proof Choose $s,t, F$ such that $H,J$ are both induced subgraphs of $F^s_s$ for some forest $F$. If $\{F^s_t, \overline{F^s_t}\}$
	has the Erd\H{o}s-Hajnal property, then so does $\{H,J\}$, so we may assume that $H=J=F^s_t$. 
	Let $d,\delta, \eta$ satisfy \ref{gethandset}. 
	By a theorem of R\"odl~\cite{rodl}, there exists $\zeta>0$ such that for every $H$-free graph $G$, there exists $S\subseteq V(G)$
	with $|S|\ge \zeta |G|$ such that one of $G[S], \overline{G}[S]$ is  $\eta/2$-sparse.
	Choose $\rho>0$ such that
	$\zeta\ge \eta^\rho$ and $\rho\ge |H|^2$.
	Choose $\alpha=1/128$, $\beta=2d+\rho+1$, $\gamma=1/12$, and choose $\delta'$ to satisfy \ref{useanticom}.
	Choose $c$ with $0<c\le \delta'$ such that $4(\rho+1)c\le 1$; we will show that 
	every $c$-critical graph contains one of $H,\overline H$ as an induced subgraph, and hence the theorem holds. 
	
	Suppose that $G$ is a $c$-critical graph that is 
	$\{H,\overline{H}\}$-free.  It follows that $|G|^c\ge 2$, and so $|G|\ge 2^{1/c}\ge 2^{4(\rho+1)}$.
	By R\"odl's theorem and replacing $G$ by its complement if necessary, we may assume that there exists 
	$S\subseteq V(G)$
	with $|S|\ge \zeta |G|$ such that $G[S]$ is $\eta$-sparse. Consequently there exists $y$ with $0<y\le \eta$ 
	such that there is a subset $S\subseteq V(G)$ with $|S|\ge y^\rho |G|$, such that $G[S]$ is $y/2$-sparse. Call such a value of $y$ {\em good}.
	Suppose there is a good value of $y$ with $|G|^{-2c}/2\le y\le |G|^{-c}$, and let $S$ be the corresponding subset. 
	Since 
	$$y|S|>y^{\rho+1}|G|\ge |G|^{1-2c(\rho+1)}2^{-\rho-1}\ge |G|^{1/2}|G|^{-1/4}\ge 2$$
	and $G[S]$
	has maximum degree at most $y|S|/2$, $S$ includes a stable set of size at least $|S|/(y|S|/2+1)\ge 1/y$ (since 
	$y|S|\ge 2$).
	But $1/y\ge |G|^c$, contradicting that $G$ is $c$-critical.
	
	On the other hand, $\eta$ is good, and $\eta>|G|^{-2c}$.
	%$\eta>|G|^{-2c}$
	Choose a good value of $y$ with $|G|^{-2c}\le y\le \eta$ and minimal with this property. It follows that $y>|G|^{-c}$, and
	$y^2/2$ is not good.
	
	From \ref{gethandset}, taking $Z=S$ and with $y$ replaced by $y/2$, we deduce that either:
	\begin{itemize}
		\item there is a subset $S'\subseteq S$ with $|S'| \ge (y/2)^{|H|^2} |S|$ that is $(y/2)^2$-sparse; or
		\item there is an equicardinal,  symmetrically $2(y/2)^{1/6}$-sparse $((y/2)^{-1/64}, (y/2)^d|S|/2)$-blockade $\mathcal{B}$ in $G$,
		and there is an $(s,t)$-handset for $\mac B$.
	\end{itemize}
	Suppose the first holds. Then $|S'|\ge (y/2)^{|H|^2} |S|\ge (y/2)^{|H|^2}y^\rho |G|\ge (y^2/2)^\rho |G|$, and so $y^2/2$ is good, a contradiction.
	Thus the second holds. Since $(y/2)^{-1/64}\ge y^{-1/128}=y^{-\alpha}$, and $(y/2)^d|S|/2\ge y^{2d+\rho}|G|/2\ge y^{\beta}|G|$, and 
	$2(y/2)^{1/6}\le y^{\gamma}$, there is a 
	symmetrically $y^\gamma$-sparse equicardinal $(y^{-\alpha}, y^\beta|G|)$-blockade $\mathcal{B}$ in $G$,
	and there is an $(s,t)$-handset for $\mac B$. By \ref{useanticom}, there is a $\mac B$-rainbow copy of $F$; and combining this
	with the handset gives a copy of $F^s_t$, a contradiction. This proves \ref{swiss2}.~\bbox


\begin{thebibliography}{99}
		
		\bibitem{APS} N. Alon, J. Pach and J. Solymosi, ``Ramsey-type theorems with forbidden subgraphs'', {\em Combinatorica} {\bf 21} (2001), 155--170.
		\bibitem{bonamy} M. Bonamy, N. Bousquet, and S. Thomass\'e, ``The Erd\H{o}s-Hajnal conjecture for long holes and antiholes'',
		{\em SIAM J. Discrete Math.}, {\bf 30} (2016), 1159--1164.
		\bibitem{loglog} M. Buci\'c, T. Nguyen, A. Scott, and P. Seymour, ``Induced subgraph density. I. A loglog step towards
		Erd\H{o}s-Hajnal'', {\em IMRN}, to appear, {\tt arXiv:2301.10147}.
		\bibitem{equiv} M. Buci\'c, J. Fox, H. T. Pham, ``Equivalence between Erd\H os-Hajnal and polynomial R\"odl and Nikiforov conjectures", arxiv:2403.08303, 2024
		\bibitem{safra} M. Chudnovsky and M. Safra, ``The Erd\H{o}s-Hajnal conjecture for bull-free graphs'', {\em J. Combinatorial Theory, Ser. B}, {\bf 98} (2008), 1301--1310.
		\bibitem{C5} M. Chudnovsky, A. Scott, P. Seymour,  and S. Spirkl,
		``Erd\H{o}s-Hajnal for graphs with no five-hole'',
		{\em Proceedings of the London Math. Soc.}, {\bf 126} (2023), 997--1014, {\tt arXiv:2102.04994}.
		\bibitem{pure1} 
		M. Chudnovsky, A. Scott, P. Seymour and S. Spirkl,
		``Pure pairs. I. Trees and linear anticomplete pairs'',
		{\em Advances in Math.}, {\bf 375} (2020), 107396, {\tt arXiv:1809.00919}.
		\bibitem{pure2} 
		M. Chudnovsky, A. Scott, P. Seymour and S. Spirkl,
		``Pure pairs. II. Excluding all subdivisions of a graph'',
		{\em Combinatorica} {\bf 41} (2021), 279--405, {\tt arXiv:1804.01060}.
		\bibitem{davies} J. Davies, ``Pivot-minors and the Erd\H{o}s-Hajnal conjecture'', {\tt arXiv:2305.09133}.
		\bibitem{EH0} P. Erd\H{o}s and A. Hajnal, ``On spanned subgraphs of graphs'',
		{\em Graphentheorie und Ihre Anwendungen} (Oberhof, 1977), 
		https://old.renyi.hu/\raisebox{-1ex}{\textasciitilde}p\_erdos/1977-19.pdf.
		\bibitem{EH}  P. Erd\H{o}s and A. Hajnal, ``Ramsey-type theorems'',
		{\em  Discrete Applied Mathematics} {\bf 25} (1989), 37--52.
		\bibitem{foxsudakov} J. Fox and B. Sudakov, ``Induced Ramsey-type theorems'', {\em Advances in Mathematics} {\bf 219} (2008), 
		1771--1800.
		\bibitem{density2} J. Fox, T. Nguyen, A. Scott and P. Seymour,
		``Induced subgraph density. II. Sparse and dense sets in cographs'', submitted for publication, {\tt arXiv:2307.00801}.
		\bibitem{density4} T. Nguyen, A. Scott and P. Seymour, ``Induced subgraph density. IV. New graphs with the Erd\H{o}s-Hajnal property'', submitted for publication, {\tt arXiv:2307.06455}.
		\bibitem{density5} T. Nguyen, A. Scott, and P. Seymour, ``Induced subgraph density. V. All paths approach Erd\H os–Hajnal'', submitted for publication,
		{\tt arXiv:2307.15032}.
		\bibitem{density6}  T. Nguyen, A. Scott, and P. Seymour, ``Induced subgraph density. VI. Bounded VC-dimension", submitted for publication, {\tt arXiv:2312.15572}.
		\bibitem{density7} T. Nguyen, A. Scott and P. Seymour, ``Induced subgraph density. VII. The five-vertex path'', submitted for publication, {\tt arXiv:2312.15333}.
		\bibitem{tungthesis} T. Nguyen, {\em Induced Subgraph Density}, Ph.D. thesis, Princeton University, May 2025, in preparation.
		\bibitem{rodl}  V. R\"odl, ``On universality of graphs with uniformly distributed edges'',
		{\em Discrete Mathematics} {\bf 59} (1986), 125--134.
	\end{thebibliography}
\end{document}